\documentclass{amsart}\usepackage{amsmath}

\usepackage{amssymb}
\usepackage{amscd}
\usepackage{thmdefs}

\input{tcilatex}
\theoremstyle{definition}
\theoremstyle{remark}
\numberwithin{equation}{section}

\input tcilatex

\begin{document}
\title[$L(F_{1})$-valued Generating Operator of $L(F_{2})$]{Trivial $L(F_{1})
$-Valued Moment Series of the Generating Operator of $L(F_{2})$}
\author{Ilwoo Cho}
\address{Dep. of Math, Univ. of Iowa, Iowa City, IA, U. S. A}
\email{ilcho@math.uiowa.edu}
\keywords{Free Group Factors, Amalgamated Moment Series.}
\maketitle

\begin{abstract}
In this paper, we will consider the generating operator $x=a+b+a^{-1}+b^{-1}$
of the free group factor $L(F_{2}),$ where $F_{2}=\,<a,b>$ is the free group
with two generators $a$ and $b.$ Let $F_{1}=\,<h>$ be a free group with one
generator $h=aba^{-1}b^{-1}$ which is group isomorphic to the integers $%
\mathbb{Z}.$ Then we can construct the free group factor $L(F_{1})$ and the
conditional expectation $E:L(F_{2})\rightarrow L(F_{1}),$ defined by $%
E\left( \underset{g\in F_{2}}{\sum }\alpha _{g}g\right) =\underset{k\in F_{1}%
}{\sum }\alpha _{k}k,$ for all $\underset{g\in F_{2}}{\sum }\alpha _{g}g$ in 
$L(F_{2}).$ Then $\left( L(F_{2}),E\right) $ is the $W^{\ast }$-probability
space with amalgamation over $L(F_{1}).$ In this paper, we will compute the 
\textbf{trivial} $L(F_{1})$-valued moment series of the generating operator $%
a+b+a^{-1}+b^{-1}$ of $L(F_{2}),$ over $L(F_{1}).$ This computation is the
good example for studying the operator-valued distribution, since the
operator-valued moment series of random variables contain algebraic and
combinatorial free probability information about opeartor-valued
distribution.
\end{abstract}

\strut 

From mid 1980's, Free Probability Theory has been developed. Here, the
classical concept of Independence in Probability theory is replaced by a
noncommutative analogue called Freeness (See [9]). There are two approaches
to study Free Probability Theory. One of them is the original analytic
approach of Voiculescu and the other one is the combinatorial approach of
Speicher and Nica (See [1], [2] and [3]). \medskip Speicher defined the free
cumulants which are the main objects in Combinatorial approach of Free
Probability Theory. The free cumulants of random variables are gotten from
the free moments of random variables via M\"{o}bius inversion. But in this
paper, we will concentrate only on computing the free moments of random
variables. And he and Nica developed free probability theory by using
Combinatorics and Lattice theory on collections of noncrossing partitions
(See [3]). Also, Speicher considered the operator-valued free probability
theory, which is also defined and observed originally by Voiculescu (See
[1]). In this paper, we will observe the important example of such
operator-valued free probability.

\strut

Let $F_{N}$ be a free group with $N$-generators and let $L(F_{N})$ be the
free group factor defined by

\strut

\begin{center}
$L(F_{N})=\overline{\mathbb{C}[F_{N}]}^{w}.$
\end{center}

\strut

In this paper, by using so-called the recurrence diagram found in [13] and
[14], we will compute the \textbf{trivial} $L(F_{1})$-valued moment series
of the ganarating operator $G$ of the free group factor $L(F_{2}),$ defined
by

\strut

\begin{center}
$G=a+b+a^{-1}+b^{-1}\in L(F_{2}),$
\end{center}

\strut

over the free group factor $L(F_{1}),$ where $F_{2}=\,<a,b>$ and $%
F_{1}=\,<aba^{-1}b^{-1}>.$ Throughout this paper, we will fix $a$ and $b$ as
the generators of the free group $F_{2}$ and we will also fix $%
h=aba^{-1}b^{-1}$ as the generator of the free group $F_{1}$ which is group
isomorphic to the integers $\mathbb{Z}.$ Let $x$ be an operator in $L(F_{2}).
$ Then there exists the Fourier expansion of $x,$

\strut

\begin{center}
$x=\underset{g\in F_{2}}{\sum }\alpha _{g}u_{g},$ \ with $\alpha _{g}\in 
\mathbb{C},$ for all $g\in F_{2}.$
\end{center}

\strut

We can regard all $g\in F_{2}$ as unitaries $u_{g}$ in $L(F_{2}).$ For the
convenience, we will denote these unitaries $u_{g}$ just by $g.$ With this
notation, it is easy to check that

$\strut $

\begin{center}
$g^{*}=u_{g}^{*}=u_{g}^{-1}=u_{g^{-1}}=g^{-1}$ in $L(F_{2}),$
\end{center}

\strut

where $g^{-1}$ is the group inverse of $g$ in $F_{2}.$ We can define the
conditional expectation $E:L(F_{2})\rightarrow L(F_{1})$ by

\strut

\begin{center}
$E\left( \underset{g\in F_{2}}{\sum }\alpha _{g}g\right) =\underset{k\in
F_{1}}{\sum }\alpha _{k}k.$
\end{center}

\strut

Then we have the $W^{*}$-probability space $\left( L(F_{2}),E\right) $ with
amalgamation over $L(F_{1}).$ Let $G$ be the generating operator $%
a+b+a^{-1}+b^{-1}$ in $L(F_{2}).$ It is easy to see that the first, second
and third trivial $L(F_{1})$-moments of $G$ vanish. i.e,

\strut

\begin{center}
$E\left( G^{k}\right) =0_{L(F_{1})},$ for $k=1,2,3,$
\end{center}

\strut

since $G^{k}$ does not contain the $h^{n}$-term, for $k=1,2,3$ and for $n\in 
\mathbb{Z},$ where $h=aba^{-1}b^{-1}$ and $h^{-1}=bab^{-1}a^{-1}.$ However,
fourth trivial $L(F_{1})$-moment $E(G^{4})$ of $G$ contains the $h$-term and
the $h^{-1}$-term. So, finding the trivial $L(F_{1})$-moments of $G$ is to
find the $h^{k}$-terms of $G^{n},$ for all $k\in \mathbb{Z}$ and $n\in 
\mathbb{N}.$

\strut

The following recurrence diagram will play a key role to find such trivial $%
L(F_{1})$-valued moment series of the generating operator $G$ of $L(F_{2})$ ;

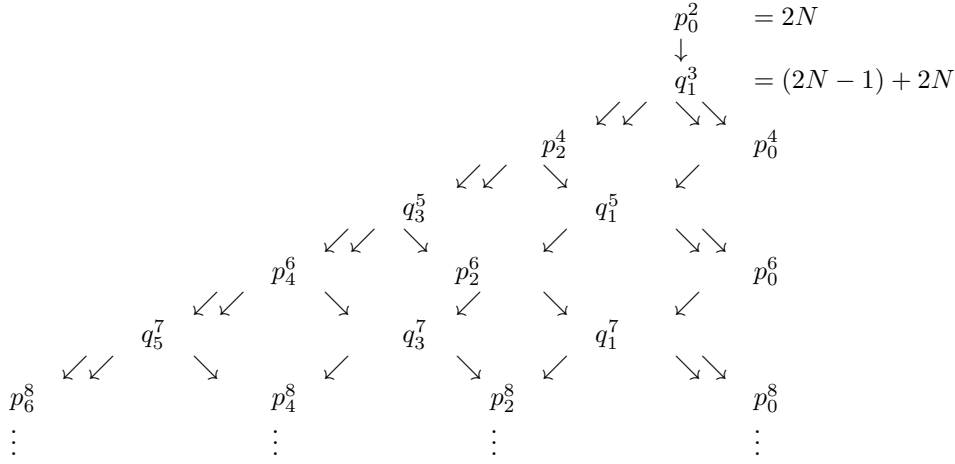
\strut

\strut

\begin{center}
$%
\begin{array}{llllllllllll}
&  &  &  &  &  &  &  &  &  & p_{0}^{2} & =2N \\ 
&  &  &  &  &  &  &  &  &  & \downarrow &  \\ 
&  &  &  &  &  &  &  &  &  & q_{1}^{3} & =(2N-1)+2N \\ 
&  &  &  &  &  &  &  &  & \swarrow \swarrow & \searrow \searrow &  \\ 
&  &  &  &  &  &  &  & p_{2}^{4} &  &  & p_{0}^{4} \\ 
&  &  &  &  &  &  & \swarrow \swarrow & \searrow &  & \swarrow &  \\ 
&  &  &  &  &  & q_{3}^{5} &  &  & q_{1}^{5} &  &  \\ 
&  &  &  &  & \swarrow \swarrow & \searrow &  & \swarrow &  & \searrow
\searrow &  \\ 
&  &  &  & p_{4}^{6} &  &  & p_{2}^{6} &  &  &  & p_{0}^{6} \\ 
&  &  & \swarrow \swarrow &  & \searrow &  & \swarrow & \searrow &  & 
\swarrow &  \\ 
&  & q_{5}^{7} &  &  &  & q_{3}^{7} &  &  & q_{1}^{7} &  &  \\ 
& \swarrow \swarrow &  & \searrow &  & \swarrow &  & \searrow & \swarrow & 
& \searrow \searrow &  \\ 
p_{6}^{8} &  &  &  & p_{4}^{8} &  &  & \text{ \ \ \ }p_{2}^{8} &  &  &  & 
p_{0}^{8} \\ 
\vdots &  &  &  & \vdots &  &  & \text{ \ \ \ }\vdots &  &  &  & \vdots%
\end{array}
$
\end{center}

\strut

where

\begin{center}
$\swarrow \swarrow $ \ : \ $(2N-1)+[$former term$]$

$\searrow $ \ \ \ \ : \ \ $(2N-1)\cdot [$former term$]$

$\swarrow $ \ \ \ \ : \ \ $\cdot +[$former term$]$
\end{center}

and

\begin{center}
$\searrow \searrow $ \ : \ \ $(2N)\cdot [$former term$].$
\end{center}

\strut

In this paper, we obtain good applications about the above recurrence
diagram. We will re-compute the moment series of the generating operators of
the free group factor $L(F_{N}),$ for all $N\in \mathbb{N},$ by using the
above recurrence diagram. This would be the one application of this
recurrence diagram (See Chapter 1). When $N=2,$ we can apply this recurrence
diagram to compute the trivial $L(F_{1})$-valued moment series of the
generating operator of $L(F_{2})$ (See Chapter 2). Remark that to study
(scalar-valued or operator-valued) moment series of elements in an operator
algebra is to study (scalar-valued or operator-valued) free distributions of
elements in that operator algebra. So, the computations in this paper about
generating operators contain the free probability information about free
distribution of those generating operators. And the free probability
information is determined by the above recurrence diagram.

\strut

In Chapter 1, we will re-compute the (scalar-valued) moment series of the
generating operator of $L(F_{N}),$ by using the recurence diagram found in
[13] and [14]. The moment series of the generating operator of $L(F_{N})$ is
already known, but here we will compute it again, by using the above
recurrence diagram. In Chapter 2, by using the reccurence diagram when $N=2$%
, we will compute the trivial $L(F_{1})$-valued moment series of the
generating operator $G=a+b+a^{-1}+b^{-1}$ of $L(F_{2}).$ Remark that the
moment series in Chapter 1 is a scalar-valued moment series and the trivial $%
L(F_{1})$-valued moment series in Chapter 2 is a operator-valued ($L(F_{1})$%
-valued) moment series.

\strut

\strut

\strut

\section{Moment Series of the Generating Operator of $L(F_{N})$}

\strut

\strut

Let $A$ be a von Neumann algebra and let $\tau :A\rightarrow \mathbb{C}$ be
the normalized faithful trace. Then we call the algebraic pair $(A,\tau ),$
the $W^{*}$-probability space and we call elements in $(A,\tau ),$ random
variables. Define the collection $\Theta _{s},$ consists of all formal
series without the constant terms in noncommutative indeterminants $%
z_{1},...,z_{s}$ ($s\in \mathbb{N}$). Then we can regard the moment series
of random variables as elements of $\Theta _{s}.$

\strut

\begin{definition}
Let $(A,\tau )$ be a $W^{*}$-probability space with its normalized faithful
trace $\tau $ and let $a\in \left( A,\tau \right) $ be a random variale. The
moment series of $a$ is defined by the formal series in $\Theta _{1},$

\strut

$\ \ \ \ \ \ \ \ \ M_{a}(z)=\sum_{n=1}^{\infty }\tau (a^{n})\,z^{n}$ .

\strut

The coefficients $\tau (a^{n})$ are called the $n$-th moments of $a,$ for
all $n\in \mathbb{N}.$
\end{definition}

\strut \strut

Let $H$ be a group and let $L(H)$ be a group von Neumann algebra. i.e,

\strut

\begin{center}
$L(H)=\overline{\mathbb{C}[H]}^{w}.$
\end{center}

\strut

Precisely, we can regard $L(H)$ as a weak-closure of the group algebra
generated by $H$ and hence

\strut

\begin{center}
$L(H)=\overline{\{\underset{g\in H}{\sum }t_{g}g:g\in H\}}^{w}.$
\end{center}

\strut

It is well known that $L(H)$ is a factor if and only if the given group $H$
is icc. (Since the free group $F_{N}$ with $N$-generators is icc, the von
Neumann group algebra $L(F_{N})$ is a factor and it is called the free group
factor.)

\strut

Now, define the canonical trace $\tau :L(H)\rightarrow \mathbb{C}$ by

\strut

\begin{center}
$\tau \left( \underset{g\in H}{\sum }t_{g}g\right) =t_{e_{H}},$ \ for all \ $%
\underset{g\in H}{\sum }t_{g}g\in L(H),$
\end{center}

\strut

where $e_{H}$ is the identity of the group $H.$ It is easy to check that the
trace $\tau $ is normalized and faithful. So, the algebraic pair $\left(
L(H),\tau \right) $ is a $W^{*}$-pobability space$.$ \strut Assume that the
group $H$ has its generators $\{g_{j}\,:\,j\in I\}.$ We say that the operator

\strut

\begin{center}
$G=\underset{j\in I}{\sum }g_{j}+\underset{j\in I}{\sum }g_{j}^{-1},$
\end{center}

\strut

the generating operator of $L(H).$ For instance, if we have a free group $%
F_{N}=\,<g_{1},...,g_{N}>,$ then the generating operator of the free group
factor $L(F_{N})$ is

\strut

\begin{center}
$g_{1}+...+g_{N}+g_{1}^{-1}+...+g_{N}^{-1}.$
\end{center}

\strut

Rest of this chapter, we will consider the moment series and the R-transform
of the generating operator $G$ of $L(F_{N}).$

\strut

From now, fix $n\in \mathbb{N}.$ And we will denote free group factor $%
L(F_{N})$ by $A$. i.e

\strut

\begin{center}
$A=\overline{\{\underset{g\in F_{N}}{\sum }t_{g}g:t_{g}\in \mathbb{C}\}}^{w}.
$
\end{center}

\strut

Recall that there is the canonical trace $\tau :A\rightarrow \mathbb{C}$
defined by

\strut

\begin{center}
$\tau \left( \underset{g\in F_{N}}{\sum }t_{g}g\right) =t_{e},$
\end{center}

\strut

where $e\in F_{N}$ is the identity of $F_{N}$ and hence $e\in L(F_{N})$ is
the unity $1_{L(F_{N})}.$ The algebraic pair $\left( L(F_{N}),\tau \right) $
is a $W^{*}$-probability space. Let $G$ be the generating operator of $%
L(F_{N}).$ i.e

\strut

\begin{center}
$G=g_{1}+...+g_{N}+g_{1}^{-1}+...+g_{N}^{-1},$
\end{center}

\strut

where \ $F_{N}=\,<g_{1},...,g_{N}>.$ It is well-known that if we denote the
sum of all words with length $n$ in $\{g_{1}$ $,$ $g_{1}^{-1}$ $,$ $...,$ $%
g_{N},$ $g_{N}^{-1}\}$ by

$\strut $

\begin{center}
$X_{n}=\underset{\left| w\right| =n}{\sum }w\in A$, for all $n\in \mathbb{N},
$
\end{center}

\strut

then

\strut

(1.1) $\ \ \ \ \ \ \ \ \ X_{1}X_{1}=X_{2}+2N\cdot e$ \ \ \ \ ($n=1$)

\strut

and

\strut

(1.2) $\ \ \ \ \ \ X_{1}X_{n}=X_{n+1}+(2N-1)X_{n-1}$ $\ \ \ \ (n\geq 2)$

\strut

(See [15]). In our case, we can regard our generating operator $G$ as $X_{1}$
in $A,$ by the very definition of $G.$

\strut

By using the relation (1.1) and (1.2), we can express $G^{n}$ in terms of $%
X_{k}$'s ; For example, $G=X_{1},$

\strut

$G^{2}=X_{1}X_{1}=X_{2}+2N\cdot e,$

\strut

$G^{3}=X_{1}\cdot X_{1}^{2}=X_{1}\left( X_{2}+(2N)e\right)
=X_{1}X_{2}+(2N)X_{1}$

$\ \ \ \ =X_{3}+(2N-1)X_{1}+(2N)X_{1}=X_{3}+\left( (2N-1)+2N\right) X_{1},$

\strut

$G^{4}=X_{4}+\left( (2N-1)+(2N-1)+2N\right) X_{2}+(2N)\left(
(2N-1)+(2N)\right) e,$

\strut

$G^{5}=X_{5}+\left( (2N-1)+(2N-1)+(2N-1)+2N\right) X_{3}$

\strut

$\ \ \ \ \ \ \ \ \ +\left( (2N-1)\left( (2N-1)+(2N-1)+(2N)\right)
+(2N)\left( (2N-1)+(2N)\right) \right) X_{1},$

\strut

$G^{6}=X_{6}+\left( (2N-1)+(2N-1)+(2N-1)+(2N-1)+2N\right) X_{4}$

\strut

$\,\,\,\ \ \ \ \ \ \ \ \ +\{(2N-1)\left( (2N-1)+(2N-1)+(2N-1)+(2N)\right) $

$\ \ \ \ \ \ \ \ \ \ \ \ \ \ \ \ \ \ \ \ \ \ \ \ \ \ \ \ \ \ \ \
+(2N-1)\left( (2N-1)+(2N-1)+(2N)\right) $

$\ \ \ \ \ \ \ \ \ \ \ \ \ \ \ \ \ \ \ \ \ \ \ \ \ \ \ \ \ \ \ \ \ \ \ \ \ \
\ \ \ \ \ \ \ \ \ \ \ \ \ \ +(2N-1)((2N-1)+(2N))\}X_{2}$

\strut

$\ \ \ \ \ \ \ \ \ \ +(2N)\left( (2N-1)\left( (2N-1)+(2N-1)+(2N))\right)
+(2N)((2N-1)+(2N))\right) e,$

etc.

\strut

So, we can find a recurrence relation to get $G^{n}$ ($n\in \mathbb{N}$)
with respect to $X_{k}$'s ($k\leq n$). Inductively, $G^{2k-1}$ and $G^{2k}$
have their representations in terms of $X_{j}$'s as follows ;

\strut

\begin{center}
$%
G^{2k-1}=X_{1}^{2k-1}=X_{2k-1}+q_{2k-3}^{2k-1}X_{2k-3}+q_{2k-5}^{2k-1}X_{2k-5}+...+q_{3}^{2k-1}X_{3}+q_{1}^{2k-1}X_{1} 
$
\end{center}

\strut \strut

and

\begin{center}
$%
G^{2k}=X_{1}^{2k}=X_{2k}+p_{2k-2}^{2k}X_{2k-2}+p_{2k-4}^{2k}X_{2k-4}+...+p_{2}^{2k}X_{2}+p_{0}^{2k}e, 
$
\end{center}

\strut

where $k\geq 2.$ Also, we have the following recurrence relation ;

\strut \strut

\begin{proposition}
Let's fix $k\in \mathbb{N}\,\setminus \,\{1\}.$ Let $q_{i}^{2k-1}$ and $%
p_{j}^{2k}$ ($i=1,3,5,...,2k-1,....$ and $j=0,2,4,...,2k,...$) be given as
before. If $p_{0}^{2}=2N$ and $q_{1}^{3}=(2N-1)+(2N)^{2},$ then we have the
following recurrence relations ;

\strut

(1) Let

$\ \ \ \ \ \ \ \ \
G^{2k-1}=X_{2k-1}+q_{2k-3}^{2k-1}X_{2k-3}+...+q_{3}^{2k-1}X_{3}+q_{1}^{2k-1}X_{1}. 
$

Then

\strut

$\ G^{2k}=X_{2k}+\left( (2N-1)+q_{2k-3}^{2k-1}\right) X_{2k-2}+\left(
(2N-1)q_{2k-3}^{2k-1}+q_{2k-5}^{2k-1}\right) X_{2k-4}$

\strut

$\ \ \ \ \ \ \ \ \ \ \ \ \ \ \ \ \ +\left(
(2N-1)q_{2k-5}^{2k-1}+q_{2k-7}^{2k-1}\right) X_{2k-6}+$

\strut

$\ \ \ \ \ \ \ \ \ \ \ \ \ \ \ \ \ +...+\left(
(2N-1)q_{3}^{2k-1}+q_{1}^{2k-1}\right) X_{2}+(2N)q_{1}^{2k-1}e.$

\strut i.e,

\strut

$\ \ \ \ \ \ \ p_{2k-2}^{2k}=(2N-1)+q_{2k-3}^{2k-1},$ $\ $

$\ \ \ \ \ \ \ p_{2k-4}^{2k}=(2N-1)q_{2k-3}^{2k-1}+q_{2k-5}^{2k-1},$

.....$...,$

\ \ \ \ \ \ \ $p_{2}^{2k}=(2N-1)q_{3}^{2k-1}+q_{1}^{2k-1}$

and

$\ \ \ \ \ \ \ \ p_{0}^{2k}=(2N)q_{1}^{2k-1}.$

\strut

(2) Let

$\ \ \ \ \ \ \ \ \ \ \ \
G^{2k}=X_{2k}+p_{2k-2}^{2k}X_{2k-2}+...+p_{2}^{2k}X_{2}+p_{0}^{2k}e.$

Then

\strut

$\ G^{2k+1}=X_{2k+1}+\left( (2N-1)+p_{2k-2}^{2k}\right) X_{2k-1}+\left(
(2N-1)p_{2k-2}^{2k}+p_{2k-4}^{2k}\right) X_{2k-3}$

\strut

$\ \ \ \ \ \ \ \ \ \ \ \ \ \ \ \ \ \ \ \ \ +\left(
(2N-1)p_{2k-4}^{2k}+p_{2k-6}^{2k}\right) X_{2k-5}+$

\strut

$\ \ \ \ \ \ \ \ \ \ \ \ \ \ \ \ \ \ \ \ \ +...+\left(
(2N-1)p_{4}^{2k}+p_{2}^{2k}\right) X_{3}+\left(
(2N-1)p_{2}^{2k}+p_{0}^{2k}\right) X_{1}.$

i.e,

\strut

$\ \ \ \ \ \ \ q_{2k-1}^{2k+1}=(2N-1)+p_{2k-2}^{2k},$ \ 

$\ \ \ \ \ \ \ q_{2k-3}^{2k+1}=(2N-1)p_{2k-2}^{2k}+p_{2k-4}^{2k},$

...$...,$ \ 

$\ \ \ \ \ \ \ q_{3}^{2k+1}=(2N-1)p_{4}^{2k}+p_{2}^{2k}$

and

$\ \ \ \ \ \ \ q_{1}^{2k+1}=(2N-1)p_{2}^{2k}+p_{0}^{2k}.$ \ \ \ 

$\square $
\end{proposition}

\strut

\strut

\begin{example}
Suppose that $N=2.$ and let $p_{0}^{2}=4$ and $q_{1}^{3}=3+p_{0}^{2}=3+4=7.$
Put

\strut

$\ \ \ \ \ \ \
G^{8}=X_{8}+p_{6}^{8}X_{6}+p_{4}^{8}X_{4}+p_{2}^{8}X_{4}+p_{0}^{8}e.$

\strut

Then, by the previous proposition, we have that

\strut

$\ p_{6}^{8}=3+q_{5}^{7},$ \ \ $p_{4}^{8}=3q_{5}^{7}+q_{3}^{7},$ \ \ $%
p_{2}^{8}=3q_{3}^{7}+q_{1}^{7}$ \ \ and \ $p_{0}^{8}=4q_{1}^{7}.$

\strut

Similarly, by the previous proposition,

\strut

\ \ \ \ $\ q_{5}^{7}=3+p_{4}^{6},$ \ \ \ $q_{3}^{7}=3p_{4}^{6}+p_{2}^{6}$ \
\ \ \ and \ \ $\ q_{1}^{7}=3p_{2}^{6}+p_{0}^{6},$

\strut

$\ \ \ \ \ \ p_{4}^{6}=3+q_{3}^{5},$ \ \ \ $p_{2}^{6}=3q_{3}^{5}+q_{1}^{5}$
\ \ \ \ \ and \ \ \ \ $p_{0}^{6}=4q_{1}^{5},$

\strut

$\ \ \ \ \ \ q_{3}^{5}=3+p_{2}^{4}$ \ \ \ \ \ \ \ and \ \ \ \ \ \ $%
q_{1}^{5}=3p_{2}^{4}+p_{2}^{4},$

\strut

$\ \ \ \ \ \ p_{2}^{4}=3+q_{1}^{3}$ \ \ \ \ \ \ \ \ \ and \ \ \ \ \ \ \ \ \ $%
p_{0}^{4}=4q_{1}^{3},$

\strut

and

$\ \ \ \ \ \ \ \ \ \ \ q_{1}^{3}=3+p_{0}^{2}=7.$

\strut

Therefore, combining all information,

\strut

$\ \ \ \ \ \ G^{8}=X_{8}+22\,X_{6}+202\,X_{4}+744\,X_{2}+1316\,e.$
\end{example}

\strut

We have the following diagram with arrows which mean that

\begin{center}
$\swarrow \swarrow $ \ : \ $(2N-1)+[$former term$]$

$\searrow $ \ \ \ \ : \ \ $(2N-1)\cdot [$former term$]$

$\swarrow $ \ \ \ \ : \ \ $\cdot +[$former term$]$
\end{center}

and

\begin{center}
$\searrow \searrow $ \ : \ \ $(2N)\cdot [$former term$].$
\end{center}

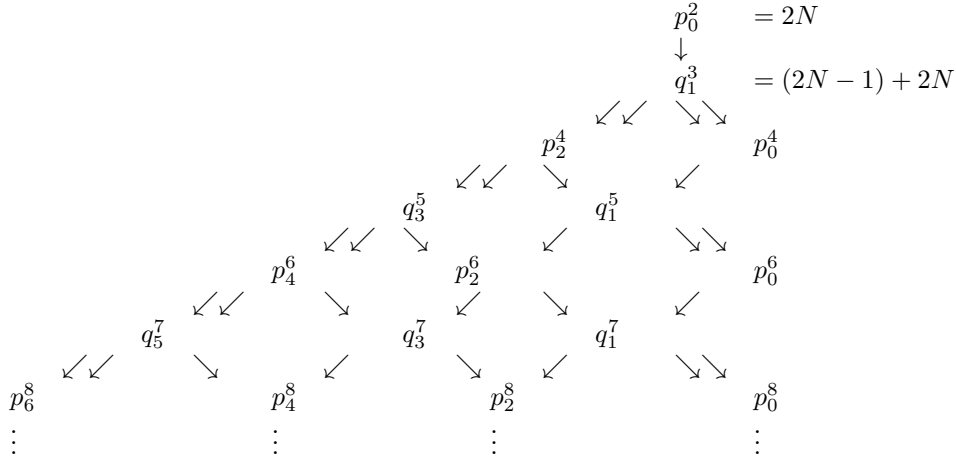
\strut

\begin{center}
$%
\begin{array}{llllllllllll}
&  &  &  &  &  &  &  &  &  & p_{0}^{2} & =2N \\ 
&  &  &  &  &  &  &  &  &  & \downarrow &  \\ 
&  &  &  &  &  &  &  &  &  & q_{1}^{3} & =(2N-1)+2N \\ 
&  &  &  &  &  &  &  &  & \swarrow \swarrow & \searrow \searrow &  \\ 
&  &  &  &  &  &  &  & p_{2}^{4} &  &  & p_{0}^{4} \\ 
&  &  &  &  &  &  & \swarrow \swarrow & \searrow &  & \swarrow &  \\ 
&  &  &  &  &  & q_{3}^{5} &  &  & q_{1}^{5} &  &  \\ 
&  &  &  &  & \swarrow \swarrow & \searrow &  & \swarrow &  & \searrow
\searrow &  \\ 
&  &  &  & p_{4}^{6} &  &  & p_{2}^{6} &  &  &  & p_{0}^{6} \\ 
&  &  & \swarrow \swarrow &  & \searrow &  & \swarrow & \searrow &  & 
\swarrow &  \\ 
&  & q_{5}^{7} &  &  &  & q_{3}^{7} &  &  & q_{1}^{7} &  &  \\ 
& \swarrow \swarrow &  & \searrow &  & \swarrow &  & \searrow & \swarrow & 
& \searrow \searrow &  \\ 
p_{6}^{8} &  &  &  & p_{4}^{8} &  &  & \text{ \ \ \ }p_{2}^{8} &  &  &  & 
p_{0}^{8} \\ 
\vdots &  &  &  & \vdots &  &  & \text{ \ \ \ }\vdots &  &  &  & \vdots%
\end{array}
$
\end{center}

\strut

\begin{quote}
\frame{\textbf{Notation}} From now, we will call the above diagram the%
\textbf{\ recurrence diagram for }$N.$ $\square $
\end{quote}

\strut

For examplet, when $N=2,$ we can compute $p_{4}^{6},$ as follows ;

\strut

\begin{center}
$p_{0}^{2}=4,$ \ \ \ \ \ $\ \ \ q_{1}^{3}=7,$
\end{center}

\strut

\begin{center}
$p_{2}^{4}=3+7=10,$ \ \ \ \ \ $p_{0}^{4}=28,$
\end{center}

\strut

\begin{center}
$q_{3}^{5}=3+10=13,$ \ $\ q_{1}^{5}=3\cdot 10+28=58.$
\end{center}

\strut

and hence $p_{4}^{6}=3+13=16.$

Recall that Nica and Speicher defined the even random variable in a $*$%
-probability space. Let $(B,\tau _{0})$ be a $*$-probability space, where $%
\tau _{0}:B\rightarrow \mathbb{C}$ is a linear functional satisfying that $%
\tau _{0}\left( b^{*}\right) =\overline{\tau _{0}(b)},$ for all $b\in B,$
and let $b\in (B,\tau _{0})$ be a random variable. We say that the random
variable $b\in (B,\tau _{0})$ is even if it is self-adjoint and it satisfies
the following moment relation ;

\strut

\begin{center}
$\tau _{0}\left( b^{n}\right) =0,$ whenever $n$ is odd.
\end{center}

\strut

By the recurrence diagram for $N$, we can get that

\strut

\begin{theorem}
Let $G\in \left( A,\tau \right) $ be the generating operator. Then the
moment series of $G$ is

\strut

$\ \ \ \ \ \ \ \ \ \tau \left( G^{n}\right) =\left\{ 
\begin{array}{lll}
0 &  & \text{if }n\text{ is odd} \\ 
&  &  \\ 
p_{0}^{n} &  & \text{if }n\text{ is even,}%
\end{array}
\right. $

\strut

for all $n\in \mathbb{N}.$
\end{theorem}

\strut

\begin{proof}
Assume that $n$ is odd. Then

\strut

$\ \ \ \ \ \
G^{n}=X_{n}+q_{n-2}^{n}X_{n-2}+...+q_{3}^{n}X_{3}+q_{1}^{n}X_{1}.$

\strut

So, $G^{n}$ does not contain the $e$-terms. Therefore,

\strut

$\ \ \ \tau \left( G^{n}\right) =\tau \left(
X_{n}+q_{n-2}^{n}X_{n-2}+...+q_{3}^{n}X_{3}+q_{1}^{n}X_{1}\right) =0.$

\strut

Assume that $n$ is even. Then

\strut

\ $\ \ \ \ G^{n}=X_{n}+p_{n-2}^{n}X_{n-2}+...+p_{2}^{n}X_{2}+p_{0}^{n}e.$

\strut

So, we have that

\strut

$\ \ \ \tau (G^{n})=\tau \left(
X_{n}+p_{n-2}^{n}X_{n-2}+...+p_{2}^{n}X_{2}+p_{0}^{n}e\right) =p_{0}^{n}.$
\end{proof}

\strut

Remark that the $n$-th moments of the generating operator in $(A,\tau )$ is
totally depending on the recurrence diagram for $N$.\strut

\strut

\begin{corollary}
Let $G\in \left( A,\tau \right) $ be the generating operator. Then $G$ is
even in $\left( A,\tau \right) $. $\square $
\end{corollary}

\strut

\begin{corollary}
Let $G\in (A,\tau )$ be the generating operator. Then the operator $G$ has
its moment series,

\strut

$\ \ \ \ \ \ \ \ \ \ \ M_{G}(z)=\sum_{n=1}^{\infty }p_{0}^{2n}\,z^{2n}\in
\Theta _{1}.$

$\square $
\end{corollary}

\strut \ \strut

\strut

\strut

\section{Trivial $L(F_{1})$-valued Moment Series of the Generating Operator
of $L(F_{2})$}

\strut

\strut

\strut

Let $M_{0}\subset M$ be von Neumann algebras with $1_{M_{0}}=1_{M}$ and let $%
\varphi :M\rightarrow M_{0}$ be the conditional expectation satisfying that

\strut

\begin{center}
$\varphi (m_{0})=m_{0},$ for all $m_{0}\in M,$
\end{center}

\strut

\begin{center}
$\varphi (m_{0}mm_{0}^{\prime })=m_{0}\cdot \varphi (m)\cdot m_{0}^{\prime
}, $
\end{center}

\strut

for all $m_{0},m_{0}^{\prime }\in M_{0},$ $m\in M,$ and

\strut

\begin{center}
$\varphi (m^{*})=\varphi (m)^{*},$ for all $m\in M.$
\end{center}

\strut

Then the algebraic pair $(M,\varphi )$ is a $W^{*}$-probability space over $%
M_{0}.$ If $m\in (M,\varphi ),$ then we will call $m$ a $M_{0}$-valued
random variable.

\strut

\begin{definition}
Let $(M,\varphi )$ be a $W^{*}$-probability space over $M_{0}$ and let $m\in
(M,\varphi )$ be a $M_{0}$-valued random variable. Define the $n$-th $M_{0}$%
-valued moment of $m$ by

\strut

$\ \ \ \ \ \ \ \ \ E\left( (m_{1}m)(m_{2}m)...(m_{n}m)\right) ,$

\strut

for all $n\in \mathbb{N},$ where $m_{1},...,m_{n}\in M_{0}$ are arbitrary.
When $m_{1}=...=m_{n}=1_{M_{0}},$ for all $n\in \mathbb{N},$ we say that the 
$M_{0}$-valued moment of $m$ is trivial. i.e, the $n$-th \textbf{trivial} $%
M_{0}$-valued moments of $m\in (M,\varphi )$ are $E(m^{n}),$ for all $n\in 
\mathbb{N}. $ We will say that the $M_{0}$-valued formal series

\strut

$\ \ \ \ \ \ \ \ \ M_{m}^{t}(z)=\sum_{n=1}^{\infty }E(m^{n})\,z^{n}\in
M_{0}[[z]]$

\strut

is the t\textbf{rivial }$M_{0}$\textbf{-valued moment series of }$m\in
(M,\varphi ),$ where $z$ is the indeterminent. (Here $M_{0}[[z]]$ is the
formal-series-ring in its indeterminent $z$)
\end{definition}

\strut \strut

In this chapter, by using the recurrence diagram for $N=2,$ we will compute
the trivial $L(F_{1})$-valued moment series of the given generating operator

\strut

\begin{center}
$G=a+b+a^{-1}+b^{-1}$
\end{center}

\strut

of the free group factor $L(F_{2}),$ where $F_{2}=\,<a,b>$ and $F_{1}=\,<h>$
are free groups, where

\strut

\begin{center}
$h=aba^{-1}b^{-1}.$
\end{center}

\strut

First, we will define the conditional expectation $E:L(F_{2})\rightarrow
L(F_{1})$ by

\strut

(2.1) $\ \ \ \ \ \ \ \ \ \ \ \ E\left( \underset{g\in F_{2}}{\sum }\alpha
_{g}g\right) =\underset{k\in F_{1}}{\sum }\alpha _{k}k,$

\strut

for all $\underset{g\in F_{2}}{\sum }\alpha _{g}g\in \left(
L(F_{2}),E\right) .$

\strut

Then we can construct the $W^{*}$-probability space $\left(
L(F_{2}),E\right) $ over its $W^{*}$-subalgebra $L(F_{1}).$ Notice that to
find the conditional expectational value of the $L(F_{1})$-valued random
variable $x\in \left( L(F_{2}),E\right) $ is to find the $h^{k}$-terms of
the $L(F_{1})$-valued random variables $x$, for $k\in \mathbb{Z},$ where $%
h=aba^{-1}b^{-1}$ and $h^{-1}=bab^{-1}a^{-1}.$

\strut

First, let us provide the recurrence diagram for $N=2$ ;

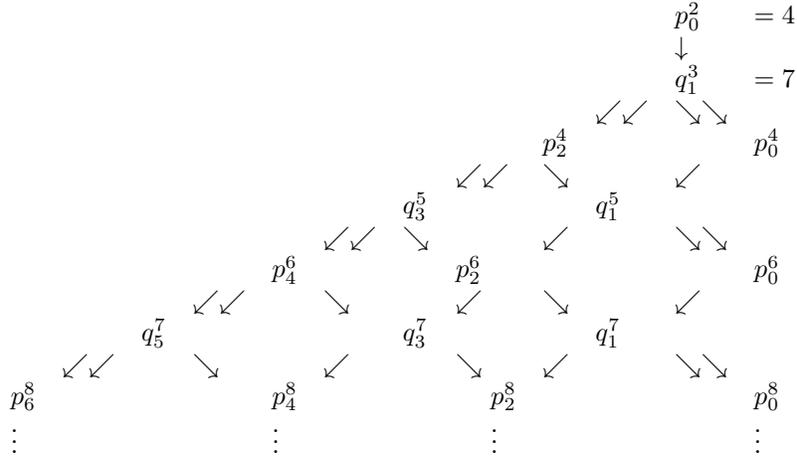
\strut

\strut

\begin{center}
$%
\begin{array}{llllllllllll}
&  &  &  &  &  &  &  &  &  & p_{0}^{2} & =4 \\ 
&  &  &  &  &  &  &  &  &  & \downarrow &  \\ 
&  &  &  &  &  &  &  &  &  & q_{1}^{3} & =7 \\ 
&  &  &  &  &  &  &  &  & \swarrow \swarrow & \searrow \searrow &  \\ 
&  &  &  &  &  &  &  & p_{2}^{4} &  &  & p_{0}^{4} \\ 
&  &  &  &  &  &  & \swarrow \swarrow & \searrow &  & \swarrow &  \\ 
&  &  &  &  &  & q_{3}^{5} &  &  & q_{1}^{5} &  &  \\ 
&  &  &  &  & \swarrow \swarrow & \searrow &  & \swarrow &  & \searrow
\searrow &  \\ 
&  &  &  & p_{4}^{6} &  &  & p_{2}^{6} &  &  &  & p_{0}^{6} \\ 
&  &  & \swarrow \swarrow &  & \searrow &  & \swarrow & \searrow &  & 
\swarrow &  \\ 
&  & q_{5}^{7} &  &  &  & q_{3}^{7} &  &  & q_{1}^{7} &  &  \\ 
& \swarrow \swarrow &  & \searrow &  & \swarrow &  & \searrow & \swarrow & 
& \searrow \searrow &  \\ 
p_{6}^{8} &  &  &  & p_{4}^{8} &  &  & \text{ \ \ \ }p_{2}^{8} &  &  &  & 
p_{0}^{8} \\ 
\vdots &  &  &  & \vdots &  &  & \text{ \ \ \ }\vdots &  &  &  & \vdots%
\end{array}
$
\end{center}

\strut

where

\begin{center}
$\swarrow \swarrow $ \ : \ $3+[$former term$]$

$\searrow $ \ \ \ \ : \ \ $3\cdot [$former term$]$

$\swarrow $ \ \ \ \ : \ \ $\cdot +[$former term$]$
\end{center}

and

\begin{center}
$\searrow \searrow $ \ : \ \ $4\cdot [$former term$].$
\end{center}

\strut

By the above recurrence diagram for $N=2,$ we have that if

$\strut $

\begin{center}
$p_{0}^{2}=4$ \ \ \ \ and $\ \ \ q_{1}^{3}=3+p_{0}^{2}=7,$
\end{center}

\strut

and if we put $X_{n}=\underset{\left| w\right| =n}{\sum }w,$ as the sum of
all words with length $n$ in $\{a,b,a^{-1},b^{-1}\},$ for $n\in \mathbb{N},$
then we have the following recurrence relations (1) and (2) ;

\strut

(1) If

\strut \strut

\begin{center}
$%
G^{2k-1}=X_{2k-1}+q_{2k-3}^{2k-1}X_{2k-3}+...+q_{3}^{2k-1}X_{3}+q_{1}^{2k-1}X_{1}. 
$
\end{center}

then

\strut

$\ \ \ \ G^{2k}=X_{2k}+\left( 3+q_{2k-3}^{2k-1}\right) X_{2k-2}+\left(
3q_{2k-3}^{2k-1}+q_{2k-5}^{2k-1}\right) X_{2k-4}$

\strut

$\ \ \ \ \ \ \ \ \ \ \ \ \ \ \ +\left(
3q_{2k-5}^{2k-1}+q_{2k-7}^{2k-1}\right) X_{2k-6}$

\strut

$\ \ \ \ \ \ \ \ \ \ \ \ \ \ \ +...+\left( 3q_{3}^{2k-1}+q_{1}^{2k-1}\right)
X_{2}+4q_{1}^{2k-1}e.$

\strut where

\strut

\begin{center}
$p_{2k-2}^{2k}=3+q_{2k-3}^{2k-1},$ $\ \ \ \ \
p_{2k-4}^{2k}=3q_{2k-3}^{2k-1}+q_{2k-5}^{2k-1},$
\end{center}

$\strut $

\begin{center}
$...,$ $p_{2}^{2k}=3q_{3}^{2k-1}+q_{1}^{2k-1}$ \ \ \ and \ \ \ $%
p_{0}^{2k}=4q_{1}^{2k-1},$
\end{center}

\strut

by the recurrence diagram.

\strut 

(2) If

\begin{center}
$G^{2k}=X_{2k}+p_{2k-2}^{2k}X_{2k-2}+...+p_{2}^{2k}X_{2}+p_{0}^{2k}e.$
\end{center}

Then

\strut

$\ \ \ \ \ \ G^{2k+1}=X_{2k+1}+\left( 3+p_{2k-2}^{2k}\right) X_{2k-1}+\left(
3p_{2k-2}^{2k}+p_{2k-4}^{2k}\right) X_{2k-3}$

\strut

$\ \ \ \ \ \ \ \ \ \ \ \ \ \ \ \ \ \ \ \ \ +\left(
3p_{2k-4}^{2k}+p_{2k-6}^{2k}\right) X_{2k-5}+$

\strut

$\ \ \ \ \ \ \ \ \ \ \ \ \ \ \ \ \ \ \ \ \ +...+\left(
3p_{4}^{2k}+p_{2}^{2k}\right) X_{3}+\left( 3p_{2}^{2k}+p_{0}^{2k}\right)
X_{1}.$

where

\strut

\begin{center}
$q_{2k-1}^{2k+1}=3+p_{2k-2}^{2k},$ \ \ \ \ \ $%
q_{2k-3}^{2k+1}=3p_{2k-2}^{2k}+p_{2k-4}^{2k},$
\end{center}

$\strut $

\begin{center}
$...,$ \ $q_{3}^{2k+1}=3p_{4}^{2k}+p_{2}^{2k}$ \ \ \ and \ \ $%
q_{1}^{2k+1}=3p_{2}^{2k}+p_{0}^{2k},$
\end{center}

\strut

by the recurrence diagram.

\strut \strut

Note that $h$ and $h^{-1}$ are words with their length 4. Therefore, $X_{4k}$
contains $h^{k}$-terms and $h^{-k}$-terms, for all $k\in \mathbb{N\cup \{}0%
\mathbb{\}}$ ! Thus we can compute the trivial $L(F_{1})$-valued moments of
the operator $G$ as follows ;

\strut

\begin{theorem}
Fix $k\in \mathbb{N}$ and Let $G\in \left( L(F_{2}),E\right) $ be the
generating operator of $L(F_{2})$. Then

\strut

(1) $\ E(G^{k})=0_{L(F_{1})},$ if $k$ is odd.

\strut

(2) $\ E\left( G^{4k}\right) =\left( h^{k}+h^{-k}\right)
+\sum_{j=1}^{k-1}p_{4k-4j}^{4k}\left( h^{k-j}+h^{-(k-j)}\right)
+p_{0}^{4k}h^{0},$

\strut

where $p_{0}^{4}=28.$

\strut

(3) If $4\nmid 2k,$ in the sense that $2k$ is not a multiple by 4, then

\strut

$\ \ \ \ \ \ E(G^{2k})=\sum_{j=1}^{k-1}p_{(2k-2)-4j}^{2k}\left( h^{\frac{k-1%
}{2}-2j}+h^{-(\frac{k-1}{2}-2j)}\right) +p_{0}^{2k}h^{0},$

\strut

where $p_{0}^{2}=4.$
\end{theorem}

\strut

\begin{proof}
(1) Suppose that $k$ is odd. Then $G^{k}$ does not have the words with
length $4p,$ for some $p\in \mathbb{N},$ by the recurrence diagram for $N=2,$
since $G^{k}$ does not have the $X_{4n}$-terms, for $n\in \mathbb{N},$ $4n<k.
$ This shows that there's no $h^{n}$-terms and $h^{-n}$-terms in $G^{k},$
where $n$ is previousely given such that $2n<k.$ Therefore, all odd trivial $%
L(F_{1})$-valued moments of $G$ vanish.

\strut

(2) By the straightforward computation using the recurrence diagram, we have
that

\strut

$\ E\left( G^{4k}\right) $

\strut

$\ \ =E\left(
X_{4k}+p_{4k-2}^{4k}X_{4k-2}+p_{4k-4}^{4k}X_{4k-4}+...+p_{4}^{4k}X_{4}+p_{2}^{4k}X_{2}+p_{0}^{4k}h^{0}\right) 
$

\strut

\ \ (2.2)\strut

\strut

$\ \ =E(X_{4k})+p_{4k-2}^{4k}E(X_{4k-2})+p_{4k-4}^{4k}E(X_{4k-4})+$

\strut

$\ \ \ \ \ \ \ \ \ \ \ \ \ \ \
...+p_{4}^{4k}E(X_{4})+p_{2}^{4k}E(X_{2})+p_{0}^{4k}h^{0}.$

\strut

Since $h^{p}$ and $h^{-p}$ terms are in $X_{4p},$ for any $p\in \mathbb{N}%
\cup \{0\},$ the formular (2.2) is

\strut

\ \ (2.3)

\strut

$\ \ \ \ \
E(X_{4k})+p_{4k-4}^{4k}E(X_{4k-4})+...+p_{4}^{4k}E(X_{4})+p_{0}^{4k}h^{0}$

\strut

$\ \ \ \ =\left( h^{k}+h^{-k}\right)
+p_{4k-4}^{4k}(h^{k-1}+h^{-(k-1)})+...+p_{4}^{4k}(h+h^{-1})+p_{0}^{4k}h^{0}.$

\strut

(3) If $4\nmid 2k,$ then $k=1,3,5,....$. If $k=1,$ then we have that ;

\strut

$\ \ \ \ \ \ \ \ \ \ \ \ \ \ \ \ \ E(G^{2})=E\left( X_{2}+4h^{0}\right)
=4h^{0}.$

\strut

If $k\neq 1$ is odd, then

\strut

$\ E(G^{2k})$

\strut

$\ \
=E(X_{2k}+p_{2k-2}^{2k}X_{2k-2}+p_{2k-4}^{2k}X_{2k-4}+p_{2k-6}^{2k}X_{2k-6}+$

\strut

$\ \ \ \ \ \ \ \ \ \ \ \ \ \ \ \ \ \ \ \ \ \ \ \ \ \ \ \
...+p_{4}^{2k}X_{4}+p_{2}^{2k}X_{2}+p_{0}^{2k}h^{0})$

\strut

$\ \
=E(X_{2k})+p_{2k-2}^{2k}E(X_{2k-2})+p_{2k-4}^{2k}E(X_{2k-4})+p_{2k-6}^{2k}E(X_{2k-6})+ 
$

\strut

$\ \ \ \ \ \ \ \ \ \ \ \ \ \ \ \ \ \ \ \ \ \ \ \ \ \ \ \
...+p_{4}^{2k}E(X_{4})+p_{2}^{2k}E(X_{2})+p_{0}^{2k}h^{0}$

\strut

$\ =0_{B}+p_{2k-2}^{2k}\left( h^{k-1}+h^{-(k-1)}\right)
+0_{B}+p_{2k-6}^{2k}\left( h^{k-3}+h^{-(k-3)}\right) +$

\strut

$\ \ \ \ \ \ \ \ \ \ \ \ \ \ \ \ \ \ \ \ \ \ \ \ \ \ \ \ \ \
...+p_{4}^{2k}(h+h^{-1})+0_{B}+p_{0}^{2k}h^{0},$

\strut

since $X_{2k-2},$ $X_{2k-6},...,X_{4}$ contain $h^{p}$-terms and $h^{-p}$%
-terms, for $p\in \mathbb{N}\cup \{0\}.$
\end{proof}

\strut

By the previous trivial $L(F_{1})$-valued moments of the generating operator 
$G$ of $L(F_{2}),$ we have the following result ;

\strut

\begin{corollary}
Let $\left( L(F_{2}),E\right) $ be the $W^{*}$-probability spcae over $%
L(F_{1})$ and let $G\in \left( L(F_{2}),E\right) $ be the generating
operator of $L(F_{2})$. Then the trivial $L(F_{1})$-valued moment series of $%
G$ is

\strut

$\ \ \ \ \ \ \ \ \ \ \ M_{G}(z)=\sum_{n=1}^{\infty }b_{2n}\,z^{n}\in
L(F_{1})[[z]],$

\strut

where

\strut

$\ \ \ \ \ \ b_{4n}=\left( h^{n}+h^{-n}\right)
+\sum_{j=1}^{n-1}p_{4n-4j}^{4n}\left( h^{n-j}+h^{-(n-j)}\right)
+p_{0}^{4n}h^{0}$

and

$\ \ \ \ \ \ \ \ b_{2k}=\sum_{j=1}^{k-1}p_{(2k-2)-4j}^{2k}\left( h^{\frac{k-1%
}{2}-2j}+h^{-(\frac{k-1}{2}-2j)}\right) +p_{0}^{2k}h^{0},$

\strut

where $4\nmid 2k,$ for all $n,k\in \mathbb{N}.$ $\square $
\end{corollary}

\strut 

\begin{remark}
Suppose we have the free group factor $L(F_{N}),$ where $N\in \mathbb{N}.$
Then we can extend the above result for the general $N.$ i.e, we can take $%
F_{N}=~<g_{1},...,g_{N}>$ \ and \ $F_{1}=~<h>,$ where

\strut 

$\ \ \ \ \ \ \ \ \ k=g_{1}\cdot \cdot \cdot g_{N}\cdot g_{1}^{-1}\cdot \cdot
\cdot g_{N}^{-1}.$

\strut 

By defining the canonical conditional expectation $E:L(F_{N})\rightarrow
L(F_{1}),$ we can construct the $W^{\ast }$-probability space with
amalgamation over $L(F_{1})$. Similar to the case when $N=2,$ by using the
recurrence diagram for $N,$ we can compute the trivial $L(F_{1})$-valued
moments of the generating operator

\strut 

$\ \ \ \ \ \ \ G_{N}=g_{1}+...+g_{N}+g_{1}^{-1}+...+g_{N}^{-1}.$

\strut 

When $N=2,$ to find the trivial $L(F_{1})$-valued moments of $G_{2}$ is to
find the $h^{k}$-terms of $G_{2}^{n},$ for $k\in \mathbb{Z},$ $n\in \mathbb{N%
}.$ But in the general case, to find the trivial $L(F_{1})$-valued moments
of $G_{N}$ is to find the $k^{p}$-terms of $G_{N}^{q},$ for $p\in \mathbb{Z},
$ $q\in \mathbb{N}.$ So, we have to choose the $X_{(2N)p}$-terms, for all $%
p\in \mathbb{N}$, containing $k^{\pm p}$-terms !
\end{remark}

\strut

\strut

\strut

\strut \textbf{References}

\bigskip 

\strut

{\small [1] \ R. Speicher, Combinatorial Theory of the Free Product with
Amalgamation and Operator-Valued Free Probability Theory, AMS Mem, Vol 132 ,
Num 627 , (1998).}

{\small [2] \ \ A. Nica, R-transform in Free Probability, IHP course note.}

{\small [3] \ \ R. Speicher, Combinatorics of Free Probability Theory IHP
course note.}

{\small [4] \ \ A. Nica, D. Shlyakhtenko and R. Speicher, R-cyclic Families
of Matrices in Free Probability, J. of Funct Anal, 188 (2002),
227-271.\strut }

{\small [5] \ \ A. Nica and R. Speicher, R-diagonal Pair-A Common Approach
to Haar Unitaries and Circular Elements, (1995), Preprint.}

{\small [6] \ \ D. Shlyakhtenko, Some Applications of Freeness with
Amalgamation, J. Reine Angew. Math, 500 (1998), 191-212.\strut }

{\small [7] \ \ A. Nica, D. Shlyakhtenko and R. Speicher, R-diagonal
Elements and Freeness with Amalgamation, Canad. J. Math. Vol 53, Num 2,
(2001) 355-381.\strut }

{\small [8] \ \ A. Nica, R-transforms of Free Joint Distributions and
Non-crossing Partitions, J. of Func. Anal, 135 (1996), 271-296.\strut }

{\small [9] \ \ D.Voiculescu, K. Dykemma and A. Nica, Free Random Variables,
CRM Monograph Series Vol 1 (1992).\strut }

{\small [10] D. Voiculescu, Operations on Certain Non-commuting
Operator-Valued Random Variables, Ast\'{e}risque, 232 (1995), 243-275.\strut 
}

{\small [11] D. Shlyakhtenko, A-Valued Semicircular Systems, J. of Funct
Anal, 166 (1999), 1-47.\strut }

{\small [12] I. Cho, The Moment Series and The R-transform of the Generating
Operator of }$L(F_{N}),${\small \ (2003), Preprint.\strut }

{\small [13] I. Cho, The Moment Series of the Generating Operator of }$%
L(F_{2})*_{L(F_{1})}L(F_{2})${\small , (2003), Preprint}

{\small [14] I. Cho, An Example of Moment Series under the Compatibility,
(2003), Preprint}

{\small [15] F. Radulescu, Singularity of the Radial Subalgebra of }$%
L(F_{N}) ${\small \ and the Puk\'{a}nszky Invariant, Pacific J. of Math,
vol. 151, No 2 (1991)\strut , 297-306.\strut \strut }

\strut

\end{document}